\documentclass[12pt]{article}

\usepackage{amssymb}

\usepackage{alltt,times,pifont,amssymb,url,graphicx} % JRH \{jrh}
\usepackage{latexsym} % RDA

\def\Hide#1{\relax}

\DeclareSymbolFont{AMSb}{U}{msb}{m}{n}
\DeclareSymbolFontAlphabet{\mathbb}{AMSb}

\def\Func#1{{\mathsf{#1}}}

\def\Cos#1{\Func{cos}\,#1}
\def\Sin#1{\Func{sin}\,#1}

\def\Imp{\Rightarrow}

\def\And{\land}

\def\Ve{\mathbf{e}}

\def\Vp{\mathbf{p}}
\def\Vq{\mathbf{q}}

\def\Vu{\mathbf{u}}
\def\Vv{\mathbf{v}}
\def\Vw{\mathbf{w}}
\def\Vx{\mathbf{x}}

\def\VO{\mathbf{0}}

\def\Norm#1{{||}#1{||}}

\newtheorem{Theorem}{Theorem}
\newtheorem{Lemma}[Theorem]{Lemma}

\def\Proof{\par \noindent{\bf Proof: }}
\def\Done{\hfill\rule{0.5em}{0.5em}}

\newcommand{\nat}{\mbox{$\protect\mathbb N$}}

\newcommand{\num}{\mbox{$\protect\mathbb Z$}}

\newcommand{\rat}{\mbox{$\protect\mathbb Q$}}
\newcommand{\real}{\mbox{$\protect\mathbb R$}}

\newcommand{\spot}{{\cdot}}

\newcommand{\all}[1]{\forall #1 \spot\:}

\newcommand{\BEQ}{\mbox{\raise4pt\hbox{$\ulcorner$}}}

\newcommand{\EEQ}{\mbox{\raise4pt\hbox{$\urcorner$}}}

\newcommand{\BA}{\begin{array}[t]{l}}
\newcommand{\EA}{\end{array}}

\makeatletter
\def\imod#1{\allowbreak\mkern10mu({\operator@font mod}\,\,#1)}
\makeatother

\newcommand\HOLSpacing{12pt}

\newlength{\hsbw}
\newcommand{\inner}[1]{\mbox{$\left\langle #1 \right\rangle$}}

\title{Aronszajn's Criterion for Euclidean Space}
\author{R.D. Arthan\thanks{This note was inspired by joint work with
Robert M. Solovay and John Harrison on decidability for logical theories
of normed spaces. I am grateful to Bob and John for ther comments.}}

\begin{document}
\maketitle
\begin{abstract}
We give a simple proof of a characterization of euclidean
space due to Aronszajn and derive a well-known characterization due to Jordan \& von Neumann as a corollary.
\end{abstract}

A norm $\Norm{\_}$ on a vector space $V$ is {\em euclidean} if
there is an inner product $\inner{\_, \_}$ on $V$ such that $\Norm{\Vv} =
\sqrt{\inner{\Vv, \Vv}}$.  Characterizations of euclidean normed spaces
abound. Amir \cite{amir86}
surveys some 350 characterizations,
starting with a well-known classic of Jordan \& von Neumann \cite{jordan35}:
a norm is euclidean iff it satisfies the {\em parallelogram identity}:
\begin{eqnarray*}
\Norm{\Vv + \Vw} &=& \sqrt{2\Norm{\Vv}^2 + 2\Norm{\Vw}^2 - \Norm{\Vv - \Vw}^2}
\end{eqnarray*}

\begin{figure}[t]
\begin{center}
\includegraphics[angle=0,scale=1.75]{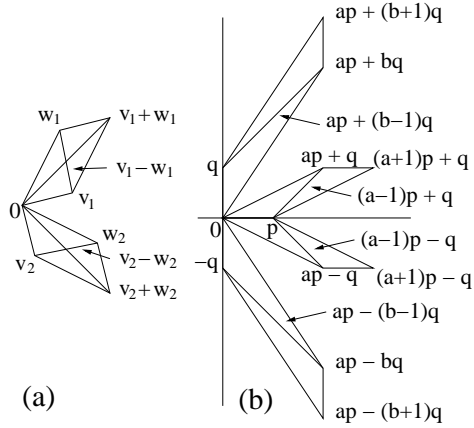}
\caption{Some parallelograms (a bird and an imp).}
\label{fig:parallelograms}
\end{center}
\end{figure}

Aronszajn proved that the algebraic details of this identity are mostly
irrelevant: if the norms of two sides and one diagonal of a parallelogram
determine the norm of the other diagonal then the norm is euclidean.
Formally, Aronszajn's criterion is the following property, as
illustrated in figure~\ref{fig:parallelograms}(a).
\[
\begin{array}{@{}l@{}l@{}}
\all{\Vv_1\;\Vw_1\;\Vv_2\;\Vw_2} &
        \Norm{\Vv_1} = \Norm{\Vv_2} \And
        \Norm{\Vw_1} = \Norm{\Vw_2} \And
        \Norm{\Vv_1 - \Vw_1} = \Norm{\Vv_2 - \Vw_2} \\
        & \quad\quad {} \Imp \Norm{\Vv_1 + \Vw_1} = \Norm{\Vv_2 + \Vw_2}.
\end{array}
\]

Aronszajn's announcement of this characterization \cite{aronszajn35} does not give a proof.
Amir's proof forms part of a long chain of interrelated results.  In this
note we give a short, self-contained proof of the theorem and derive the
Jordan-von Neumann theorem as a corollary.
We begin with a lemma showing that the Aronszajn
criterion ensures a useful supply of isometries.
Figure~\ref{fig:parallelograms}(b) illustrates the parallelograms
that feature  in the proof.
\begin{Lemma}\label{lemma:refl}
Let $V$ be a 2-dimensional normed space satisfying the Aronszajn criterion.
If $\VO \not= \Vp, \Vq \in V$ and
$\Norm{\Vp + \Vq} = \Norm{\Vp - \Vq}$, then there is a linear isometry
$\mu:V \rightarrow V$ such that $\mu(\Vp) = \Vp$ and $\mu(\Vq) = -\Vq$.
\end{Lemma}
\Proof
The equation $\Norm{\Vp - \Vq} = \Norm{\Vp + \Vq}$ cannot hold if the non-zero 
vectors $\Vp$ and $\Vq$ are linearly dependent, so $\Vp$ and $\Vq$
form a basis for $V$ and
$\mu(a\Vp + b\Vq) = a\Vp - b\Vq$ will define the required isometry provided
the following holds  for all $a, b \in \real$.
\[
\begin{array}{crcl}
\hspace{1.35in}& \Norm{a\Vp + b\Vq} &=& \Norm{a\Vp - b\Vq} \hspace{1in} \mbox{(*)}
\end{array}
\]
(*) is trivial when $a = 0$ or $b = 0$ and is true by assumption when $a = b = 1$.
The instances of the Aronszajn criterion displayed below hold in $V$ by assumption.
As the antecedent of the implication may be assumed by induction for integer
$a > 1$,
we conclude that (*) holds for $a \in \nat$ and $b = 1$.
\[
\begin{array}{crcl}
      & \Norm{a\Vp + \Vq}       &=& \Norm{a\Vp - \Vq}       \\
 \And & \Norm{\Vp}              &=& \Norm{\Vp}              \\
 \And & \Norm{(a-1)\Vp + \Vq}   &=& \Norm{(a-1)\Vp - \Vq}   \\
 \Imp & \Norm{(a+1)\Vp + \Vq}   &=& \Norm{(a+1)\Vp - \Vq}
\end{array}
\]
This gives the base case for an induction on $b$
showing that (*) holds for $a, b \in \nat$
using the following instances of the Aronszajn criterion.
\[
\begin{array}{crcl}
      & \Norm{a\Vp + b\Vq}       &=& \Norm{a\Vp - b\Vq}     \\
 \And & \Norm{\Vq}               &=& \Norm{{-}\Vq}          \\
 \And & \Norm{a\Vp + (b-1)\Vq}   &=& \Norm{a\Vp - (b-1)\Vq} \\
 \Imp & \Norm{a\Vp + (b+1)\Vq}   &=& \Norm{a\Vp - (b+1)\Vq}
\end{array}
\]

By symmetry, (*) holds for all $a, b \in \num$; using $\Norm{(j/k)\Vp +
(m/n)\Vq} = |1/(kn)|\cdot\Norm{jn\Vp + km\Vq}$, we find that (*) holds for all
$a, b \in \rat$; finally, by continuity, (*) holds for all $a, b \in \real$.
\Done

The next lemma shows that the conclusion of lemma~\ref{lemma:refl}
characterizes euclidean space amongst 2-dimensional normed spaces.
\begin{Lemma}\label{lemma:2d}
Let $V$ be a 2-dimensional normed space such that if
$\VO \not= \Vp, \Vq \in V$ and
$\Norm{\Vp + \Vq} = \Norm{\Vp - \Vq}$, then there is a linear isometry
$\mu:V \rightarrow V$ with $\mu(\Vp) = \Vp$ and $\mu(\Vq) = -\Vq$.
Then $V$ is euclidean.
\end{Lemma}
\Proof
Let $\Ve_1 \in V$ be a unit vector. As $\Vx$ traverses an arc of the
unit circle from $\Ve_1$ to $-\Ve_1$, $\Norm{\Ve_1 + \Vx} - \Norm{\Ve_1
- \Vx}$ varies continuously from $2$ to $-2$ and hence is 0 for some
  $\Vx$. Let $\Ve_2$ be such an $\Vx$.  Take euclidean coordinates with
respect to $\Ve_1$ and $\Ve_2$ and so fix a euclidean norm $\Norm{\_}_e$ on
$V$ together with the associated notions of angle, rotation and reflection
in a line.
The condition $\Norm{\Vp + \Vq} = \Norm{\Vp - \Vq}$ is satisfied both for
$\Vp = \Ve_1$, $\Vq = \Ve_2$ and for $\Vp = \Ve_1 + \Ve_2$, $\Vq = \Ve_1 -
\Ve_2$, and so, by assumption, the reflections in the $x$-axis and in the line $x
= y$ are both $V$-isometries, and hence so is their composite, a rotation
$\rho$ through a right angle.
Let $\Vu$ be any unit vector in $V$, let $\Vp$ bisect
the angle between $\Vu$ and $\Ve_1$ and let $\Vq = \rho(\Vp)$.
We then have:
\[
\Norm{\Vp + \Vq} = \Norm{\rho(\Vp + \Vq)} = \Norm{\rho(\Vp) + \rho(\Vq)} = \Norm{\Vq - \Vp} =
\Norm{\Vp - \Vq}
\]
\noindent and so our assumptions imply that the reflection $\mu$ in the
bisector of the angle between $\Vu$ and $\Ve_1$
is a $V$-isometry. Thus $\Vu$ and $\mu(\Ve_1)$ have the
same direction and the same $V$-norm, and so $\Vu = \mu(\Ve_1)$.
Hence $\Norm{\Vu}_e = \Norm{\mu(\Ve_1)}_e = \Norm{\Ve_1}_e = 1$, for every unit vector $\Vu$ of $V$, and $V$ is therefore euclidean.
\Done

To lift results from 2 dimensions to arbitrary dimensions we use the
following lemma, which Amir proves using the Jordan-von Neumann
theorem, as did Jordan \& von Neumann.  We give a geometric proof.
\begin{Lemma}\label{lemma:2d-implies-any}
If every 2-dimensional subspace of a normed space $V$ is euclidean,
then $V$ is euclidean.
\end{Lemma}
\Proof
Let $d$ be the (possibly infinite) dimension of $V$.
If $d = 0, 1$ or $2$, the theorem is trivially true, so we may
assume $d \ge 3$.
If $V$ is euclidean, then easy algebra shows that the inner product
$\inner{\_, \_}$ must be given by:
\[
\inner{\Vv, \Vw} = (\Norm{\Vv + \Vw}^2 - \Norm{\Vv}^2 - \Norm{\Vw}^2)/2.
\]
So $V$ is euclidean iff the function $\inner{\Vv, \Vw}$
defined by the above equation satisfies the axioms for an inner product.
These axioms can be expressed using just 3 vector variables,
and so they hold in $V$ iff they hold in every subspace of $V$ of
dimension at most 3. Thus we may assume $d = 3$.

Let $U$ be a 2-dimensional and so euclidean subspace of $V$;
let $\Ve_1$ and $\Ve_2$ be orthogonal unit vectors in $U$;
and let $\Ve_3$ lie in the intersection of the unit sphere $S_V$
and a supporting plane of $S_V$, say $P = \Ve_3 + U$, parallel to $U$.
For $\theta \in [0, \pi)$,
let $\Vp_{\theta}$ be the unit vector $\Cos\theta\cdot\Ve_1 + \Sin\theta\cdot\Ve_2$
in $U$ and let $W_{\theta}$ be the subspace of $V$
spanned by $\Vp_{\theta}$ and $\Ve_3$.
$W_{\theta}$ is 2-dimensional and so euclidean.
The line $P \cap W_{\theta}$ is parallel to $\Vp_{\theta}$ and
meets $S_{W_{\theta}}$ at $\Ve_3$.
Hence as $P$ supports $S_V$, $P \cap W_{\theta}$ must
support the unit circle $S_{W_{\theta}}$ of $W_{\theta}$.
It follows that $\Ve_3$ is a unit vector orthogonal to $\Vp_{\theta}$ in
$W_{\theta}$
and that
$S_{W_{\theta}} = \{ \Sin\phi\cdot\Vp_{\theta} + \Cos\phi\cdot\Ve_3 \mid \phi \in [0, 2\pi)\}$.
But $S_V$ is the union of the $S_{W_{\theta}}$, so we have:
\[
S_V = \{ \Sin\phi\cdot\Cos\theta\cdot\Ve_1 + \Sin\phi\cdot\Sin\theta\cdot\Ve_2 + \Cos\phi\cdot\Ve_3
         \mid \theta \in [0, \pi), \phi \in [0, 2\pi) \}
\]
I.e., the unit sphere of $V$ is a euclidean sphere, and $V$ must be euclidean.
\Done
\begin{Theorem}
[Aronszajn; Jordan \& von Neumann]
\label{thm:aronszajn}
\label{thm:jordan-von-neumann}
If $V$ is a normed space, then the following are equivalent:
\begin{center}
\begin{tabular}{r@{$\;$}l}
{\it(i)} & the Aronszajn criterion holds in $V$; \\
{\it(ii)} & $V$ is euclidean; \\
{\it(iii)} & the parallelogram identity holds in $V$.
\end{tabular}
\end{center}
\end{Theorem}
\Proof
[$\mbox{{\it(i)}} \Imp \mbox{{\it(ii)}}$]: this is immediate
from lemmas~\ref{lemma:refl}, \ref{lemma:2d} and~\ref{lemma:2d-implies-any};\\\relax
[$\mbox{{\it(ii)}} \Imp \mbox{{\it(iii)}}$]: easy algebra 
shows that with $\Norm{\Vv} = \sqrt{\inner{\Vv, \Vv}}$ the
parallelogram identity holds in an inner product space;\\\relax
[$\mbox{{\it(iii)}} \Imp \mbox{{\it(i)}}$]:  the parallelogram identity clearly implies the Aronszajn criterion.
\Done

\bibliographystyle{plain}
\bibliography{bookspapers}

\newpage

\end{document}